# FIBRES, CONNEXIONS ET HOMOLOGIE CYCLIQUE
## par Max KAROUBI

## 1. La notion d'espace fibré au sens de Steenrod

**1.1.** Dans l'acceptation la plus simple, un espace fibré $\xi$ de base B et de fibre F est la donnée d'une application continue (notée traditionnellement de manière verticale)

$$\begin{array}{c} E \\ \downarrow \pi \\ B \end{array}$$

telle que $\forall\, b \in B$, il existe un voisinage ouvert de b - soit $V = V(b)$ - et un homéomorphisme $\varphi : V \times F \longrightarrow \pi^{-1}(V)$, tel que le diagramme suivant commute

$$\begin{array}{ccc} V \times F & \longrightarrow & \pi^{-1}(V) \\ & \searrow \quad \swarrow & \\ & V & \end{array}$$

La flèche oblique de gauche est la première projection et celle de droite la restriction de $\pi$ à $\pi^{-1}(V)$. En particulier, si on peut choisir $V = B$, on dit que le fibré est trivial. On note E "l'espace total" du fibré qu'on identifiera par abus de langage au fibré $\xi$. Par ailleurs, si C est un sous-espace de B, la "restriction" de E à C est simplement la restriction de l'application continue $\pi$ à $\pi^{-1}(C)$ : c'est un fibré de base C.

**1.2.** Soit $(U_i)$ un recouvrement ouvert de B tel que la restriction de E à chaque $U_i$ soit un fibré trivial (le recouvrement est alors dit "trivialisant"). On a donc des homéomorphismes

$$\varphi_i : U_i \times F \longrightarrow \pi^{-1}(U_i)$$

d'où des "données de recollement"

$$(\varphi_j.\varphi_i^{-1}) : (U_i \cap U_j) \times F \longrightarrow (U_i \cap U_j) \times F$$

L'homéomorphisme $\varphi_j.\varphi_i^{-1}$ est de la forme

$$(x, e) \mapsto (x, g_{ji}(x)(e))$$

où $g_{ji} : U_i \cap U_j \longrightarrow \mathrm{Aut}(F)$ est une certaine application. Un G-fibré "repéré" ("coordinate bundle" au sens de Steenrod [S]) est la donnée d'un sous-groupe topologique G de Aut(F) et de trivialisations $(\varphi_i)$ comme ci-dessus tels que les $g_{ji}$ définissent des applications <u>continues</u> de $U_i \cap U_j$ dans $G \subset \text{Aut}(F)$. On dit alors que G est le "groupe structural" du fibré E.

**1.3**. Il convient de remarquer que cette condition de continuité des $g_{ji}$ n'est pas automatique, sauf dans des cas particuliers. Par exemple, si E est un fibré vectoriel (structure linéaire sur chaque fibre, trivialisations compatibles avec ces structures linéaires), il est naturel de prendre pour G le groupe linéaire $GL_n(\mathbf{C})$ ou $GL_n(\mathbf{R})$, suivant que l'on considère des fibrés complexes ou réels. Dans ce cas, la condition de continuité est automatique.

Il en va différemment si les fibres sont des espaces vectoriels de dimension infinie (par exemple des espaces de Banach F) : la continuité de l'application

$$(x,, e) \mapsto (x, g_{ji}(x)(e))$$

n'implique pas nécessairement celle de l'application $x \mapsto g_{ji}(x)$ de $U_i \cap U_j$ dans G, où G est le groupe banachique des automorphismes linéaires continus de F.

Pour un recouvrement ouvert trivialisant $\mathcal{U}$ donné, il est facile de voir qu'on définit bien ainsi une "catégorie" de G-fibrés. Ce qui nous intéresse surtout, ce sont les classes d'isomorphie de tels fibrés. Steenrod a montré que les "fonctions de transition" $g_{ji}$ caractérisent le G-fibré E. De manière précise, on vérifie d'abord la "condition de cocycle"

$$g_{ki} = g_{kj} \cdot g_{ji}$$

au-dessus de chaque point $x \in U_i \cap U_j \cap U_k$. Réciproquement, la donnée d'un tel cocycle permet de reconstruire le fibré E comme quotient de la réunion disjointe des $U_i \times F$ par la relation d'équivalence évidente qui identifie $(x_i, e)$ à $(x_j, g_{ji}(x)(e))$ pour $x_i \in U_i$, $x_j \in U_j$ et $x = x_i = x_j \in U_i \cap U_j$. Deux cocycles $(g_{ji})$ et $(h_{ji})$ définissent deux fibrés isomorphes s'ils sont "cohomologues" (au sens non abélien) : ceci veut dire qu'il existe des fonctions continues $\lambda_i : U_i \longrightarrow G$ telles que

$$h_{ji} = (\lambda_j)^{-1} g_{ji} \lambda_i$$

au-dessus d'un point $x \in U_i \cap U_j$. On désigne par $H^1(\mathcal{U} ; G)$ l'ensemble des classes d'équivalence de tels cocycles (ensemble de "cohomologie non abélienne").

**1.4.** Un autre recouvrement ouvert $\mathcal{V}$ de X est dit plus fin que $\mathcal{U}$ si, pour tout $V \in \mathcal{V}$, il existe $U \in \mathcal{U}$ avec $V \subset U$. Ceci permet de définir une application

$$H^1(\mathcal{U} ; G) \longrightarrow H^1(\mathcal{V} ; G)$$

de la manière suivante. Choisissons pour chaque $V_\alpha \in \mathcal{V}$ un $U_{i(\alpha)} \in \mathcal{U}$ avec $V_\alpha \subset U_{i(\alpha)}$. On définit alors une nouvelle fonction de transition $\theta_{\beta\alpha}$ comme la restriction à $V_\alpha \cap V_\beta$ de $g_{j(\beta)i(\alpha)}$. On montre que la correspondance $g_{ji} \mapsto \theta_{\beta\alpha}$ induit une application de

$H^1(\mathcal{U} ; G)$ dans $H^1(\mathcal{V} ; G)$ qui est indépendante des différents choix. Avec des précautions logiques (par exemple en considérant des recouvrements ouverts indexés par les parties de X), on peut considérer la limite inductive suivant $\mathcal{U}$ des $H^1(\mathcal{U} ; G)$, limite qu'on note $H^1(X ; G)$. Si G est discret abélien, on retrouve le premier groupe de cohomologie de Cech classique.

**1.5.** Le calcul "explicite" de cet ensemble $H^1(X ; G)$ n'est pas du tout évident, même pour des exemples aussi simples qu'une sphère X de dimension > 1. Steenrod démontre [S] que cette classification des G-fibrés repérés se ramène à un problème de topologie algébrique (à priori aussi difficile à résoudre...). De manière précise, supposons que G soit un sous-groupe fermé d'un groupe linéaire, par exemple $GL_n(\mathbf{C})$. Steenrod considère alors l'application quotient

$$E_m = GL_{n+m}(\mathbf{C})/GL_m(\mathbf{C}) \longrightarrow GL_{n+m}(\mathbf{C})/G \times GL_m(\mathbf{C}) = B_m$$

L'espace $E_m$ est une "variété de Stiefel", tandis que $B_m$ est l'analogue d'une grassmanienne (c'est exactement une grassmanienne si G est le groupe linéaire $GL_n(\mathbf{C})$). Grâce à la théorie des groupes de Lie, on peut montrer que les $E_m$ sont "universels" dans le sens suivant : considérons la limite inductive des espaces $E_m$ et $B_m$. On obtient alors une fibration $E \longrightarrow B$, où la base B est ce qu'on appelle l'espace classifiant du groupe G, noté BG. Il vérifie la propriété remarquable suivante : pour tout espace paracompact X, il existe une correspondance bijective naturelle entre $H^1(X ; G)$ et l'ensemble des classes d'homotopie d'applications continues de X dans l'espace classifiant BG.

Choisissons par exemple pour X une sphère $S^r$. On peut alors montrer directement ou par le résultat précédent que $H^1(X ; G)$ est le quotient du groupe d'homotopie $\pi_{r-1}(G)$ par l'action du groupe $\pi_0(G)$ opérant par automorphisme intérieur. Toute information sur ce groupe d'homotopie permet alors une connaissance des fibrés sur des sphères. Un exemple célèbre est celui où $G = GL_n(\mathbf{R})$ avec n > r. L'ensemble $H^1(X ; G)$ peut être alors muni d'une structure de groupe abélien périodique de période 8 par rapport à r (théorème de périodicité de Bott [B]). Voici les 8 premiers groupes

| r | 0 | 1 | 2 | 3 | 4 | 5 | 6 | 7 |
|---|---|---|---|---|---|---|---|---|
| $H^1(S^r ; G)$ | **Z** | **Z**/2 | **Z**/2 | 0 | **Z** | 0 | 0 | 0 |

On peut de même calculer $H^1(S^r ; GL_m(\mathbf{C}))$ si m > r/2. On trouve cette fois-ci des groupes périodiques de période 2 par rapport à r : **Z**, 0, **Z**, 0, etc.

## 2. Connexion et courbure sur un fibré différentiable.

**2.1.** Nous nous sommes placés jusqu'à maintenant dans un cadre "topologique" : espaces topologiques et applications continues. Cependant, en géométrie on considère plutôt le cadre "différentiable" : variétés de classe $C^\infty$ et applications différentiables de classe $C^\infty$. En suivant les mêmes idées, on peut classifier les fibrés différentiables à isomorphisme près (G étant un groupe de Lie). On trouve alors un nouvel ensemble de cohomologie non abélienne qu'il est naturel de noter $H^1_{diff}(X\,;G)$. On a évidemment une application naturelle

$$\theta : H^1_{diff}(X\,;G) \longrightarrow H^1(X\,;G)$$

En utilisant des méthodes d'approximation, il est possible de montrer que cette application $\theta$ est bijective (si X est une variété paracompacte). Steenrod considérait ce résultat avec la signification suivante : on ramène la classification différentiable des fibrés à une classification continue qui se "résout" par des méthodes homotopiques standard (cf. l'exemple ci-dessus des fibrés vectoriels sur des sphères). Dans l'esprit de Steenrod, il devait être plus difficile de construire des invariants de fibrés différentiables que de fibrés continus.

**2.2.** Nous allons maintenant inverser le point de vue de Steenrod et montrer qu'il est en fait plus facile de construire des invariants de fibrés différentiables que de fibrés continus. Pour cela, la notion de connexion est fondamentale et nous allons choisir le point de vue le plus "terre à terre" possible.

En suivant de nouveau les idées de Steenrod, et en considérant pour simplifier le groupe $G = GL_n(\mathbf{C})$, nous définissons une "<u>connexion repérée</u>" comme la donnée d'un recouvrement trivialisant $\mathcal{U}$ et pour chaque ouvert $U_i$ du recouvrement, d'une matrice $\Gamma_i$ de formes différentielles de degré 1, soit $\Gamma_i \in M_n(\Omega^1(U_i))$. Ces matrices se recollent par la formule suivante (qui a un sens au-dessus de $U_i \cap U_j$)

$$\Gamma_i = g_{ij}\,\Gamma_j\,g_{ji} + g_{ij}\,dg_{ji}$$

où les $g_{ji}$ sont les fonctions de transition du fibré. Dans cette formule, on interprète l'expression $g_{ij}\,dg_{ji}$ comme l'image réciproque de la "forme de Maurer-Cartan" $g^{-1}\,dg$ définie sur le groupe de Lie $G = GL_n(\mathbf{C})$ par l'application

$$g_{ji} : U_i \cap U_j \longrightarrow G$$

Il existe plusieurs méthodes pour construire une connexion sur un fibré. Peut être la plus simple est de considérer une partition différentiable de l'unité $(\alpha_k)$ associée au recouvrement $\mathcal{U}$ et de poser

$$\Gamma_i = \sum_k \alpha_k \, g_{ik} \, dg_{ki}$$

("barycentre" de formes de Maurer-Cartan). Il est facile de vérifier que ces matrices $\Gamma_i$ de formes différentielles se recollent comme il a été indiqué dans la formule ci-dessus.

**2.3.** Nous définissons la <u>courbure</u> d'une connexion repérée comme étant définie "localement" par des matrice $R_i$ de formes différentielles de degré 2, soit

$$R_i = d\Gamma_i + (\Gamma_i)^2$$

On peut remarquer que l'expression $(\Gamma_i)^2$ peut aussi s'écrire $\frac{1}{2}[\Gamma_i, \Gamma_i]$, la moitié du commutateur <u>gradué</u> de $\Gamma_i$ avec lui-même. La courbure (contrairement à la connexion) est un "tenseur", ce qui veut dire qu'on a une formule de changement de cartes nettement plus simple

$$R_i = (g_{ji})^{-1} R_j (g_{ji})$$

**2.4. DEFINITION/THEOREME.** *Le caractère de Chern de degré* 2p *du fibré* E *muni de la connexion* $\Gamma$ *est donné par la formule suivante*

$$Ch_p(E, \Gamma) = \frac{1}{(2i\pi)^p \, p!} \, Trace(R_i)^p$$

*C'est une forme différentielle fermée indépendante de* i *qui définit donc une classe de cohomologie de* de Rham *globale appartenant à* $H^{2p}(X ; \mathbf{C})$. *Cette classe de cohomologie est indépendante du choix de la connexion sur le fibré* E.

**2.5. Remarques**. En faisant la somme pour tous les p, on voit ainsi que le caractère de Chern définit un morphisme d'ensemble (avec $G = GL_n(\mathbf{C})$)

$$H^1(X ; G) \longrightarrow H^{pair}(X ; \mathbf{C})$$

qui permet en particulier de détecter la non trivialité de certains fibrés. Par ailleurs, on peut montrer (ce n'est pas évident) que l'élement ainsi trouvé dans $H^{pair}(X ; \mathbf{C})$ est en fait dans l'image de l'homomorphisme $H^{pair}(X ; \mathbf{Z}) \longrightarrow H^{pair}(X ; \mathbf{C})$. Le facteur de normalisation $\frac{1}{(2i\pi)^p \, p!}$ a été mis dans la formule du caractère de Chern afin qu'il en soit ainsi.

# 3. Du caractère de Chern à l'homologie cyclique.

**3.1.** Supposons pour simplifier que le fibré E du § 2 soit trivialisé au-dessus d'un nombre fini d'ouverts $U_i$ i = 1, ..., r (ceci est le cas si la variété est paracompacte avec un nombre fini de composantes connexes), la fibre étant isomorphe à $\mathbf{C}^n$. On peut alors montrer que E est facteur direct d'un fibré trivial ; en d'autres termes, il existe une application différentiable

$$Q : X \longrightarrow M_N(\mathbf{C})$$

telle que $Q^2 = Q$ (plus précisément $Q(x)^2 = Q(x)$ pour tout $x \in X$). Le fibré E n'est autre que "l'image" de Q, c'est-à-dire l'ensemble des couples (x, Q(x)v) où $x \in X$ et $v \in \mathbf{C}^N$. Une formule explicite pour Q est la suivante : si $(\alpha_i)$ est une partition de l'unité associée au recouvrement $\mathcal{U} = (U_i)$, posons $\beta_i = \alpha_i / \sqrt{\sum (\alpha_j)^2}$. Alors $(\beta_j)$ est une "partition carrée" de l'unité, c'est-à-dire qu'on a notamment $\sum (\beta_i)^2 = 1$. Considérons alors la matrice Q en r x r blocs, celui situé à la place (i, j) étant égal à $\beta_i \beta_j g_{ji}$. Notons que cette matrice est de type N x N avec N = r.n .

**3.2. THEOREME**. *La matrice* Q *est un projecteur. Son image est le fibré vectoriel* E *associé aux fonctions de transition* $(g_{ji})$.

**3.3.** En considérant la "connexion de Levi-Civita" associée à ce projecteur, on peut alors montrer que le caractère de Chern de degré 2p est défini par la formule remarquable suivante (cf. [ C] et [ K]) :

$$Ch_p(E) = \frac{1}{(2i\pi)^p p!} \text{Trace } (Q.(dQ)^{2p})$$

**3.4.** L'intérêt essentiel de cette formule est qu'elle permet d'élargir la théorie des classes caractéristiques de fibrés vectoriels à la "géométrie non commutative" [ C] . De manière précise, considérons une $\mathbf{C}$-algèbre A non nécessairement commutative (on peut remplacer $\mathbf{C}$ par un anneau commutatif quelconque contenant le corps des rationnels, modulo quelques conventions de normalisation). On suppose aussi donnée une algèbre différentielle graduée $\Omega^*(A)$ dont la partie de degré 0 est précisément l'algèbre A.

Un exemple qu'on peut avoir à l'esprit est celui où $A = M_N(C^\infty(X))$, l'algèbre des matrices N x N dont les éléments sont des fonctions $C^\infty$ sur X. Dans ce cas $\Omega^r(A)$ est le A-bimodule des matrices N x N dont les coefficients sont des formes différentielles de degré r.

En partant de l'algèbre A, il existe cependant un exemple "universel" d'algèbre différentielle graduée $\Omega^*(A)$ ; on choisit $\Omega^1(A) = \text{Ker } (A \otimes A \longrightarrow A)$, où la flèche est la multiplication. On obtient ainsi un A-bimodule dont le A-produit tensoriel n fois par lui-même permet de définir $\Omega^n(A)$. Par ailleurs, la différentielle

$$d : A \longrightarrow \Omega^1(A)$$

définie par $d(a) = 1 \otimes a - a \otimes 1$ vérifie l'identité de Leibniz (non commutative)

$$d(ab) = d(a).b + a.d(b)$$

et elle s'étend à tous les $\Omega^n(A)$ suivant une méthode standard (cf. [ K]). On désigne par $\overline{\Omega}^*(A)$ le complexe quotient de $\Omega^*(A)$ par le sous-complexe engendré (additivement) par les commutateurs gradués $\omega_q \theta_p - (-1)^{qp} \theta_p \omega_q$ où $\omega_q$ est de degré q et $\theta_p$ de degré p. La cohomologie en degré n de ce complexe quotient est par définition "l'homologie non commutative" en degré n de l'algèbre A. Cette homologie $\overline{H}_n(A)$ est très proche de l'homologie cyclique $HC_n(A)$ définie par A. Connes. Le résultat précis est le suivant (cf. [C][K]) :

**3.5. THEOREME**. *On a un isomorphisme naturel*

$$\overline{H}_n(A) \cong \text{Ker } [\overline{HC}_n(A) \xrightarrow{B} HH_{n+1}(A) ]$$

*où* $\overline{HC}_n(A)$ *est l'homologie cyclique réduite et* $HH_n(A)$ *l'homologie de Hochschild.*

**3.6. DEFINITION**. *Soit E un A-module projectif de type fini égal à l'image d'un projecteur Q dans $A^n$. Alors le caractère de Chern $Ch_p(E)$ est défini par la formule précédente*

$$Ch_p(E) = \frac{1}{(2i\pi)^p \, p!} \text{ Trace } (Q.(dQ)^{2p})$$

*C'est un élément de* $\overline{H}_{2p}(A) \subset \overline{HC}_{2p}(A)$

**3.7**. Les considérations précédentes permettent d'élargir considérablement la théorie classique des classes caractéristiques relatives à un groupe de Lie G de la manière suivante. Supposons que G soit un groupe de Lie Fréchet et soit $\rho : G \longrightarrow GL_n(A)$ une représentation différentiable, où A est une algèbre de Fréchet (un cadre plus général est possible). Si E est un G-fibré différentiable sur une variété X au sens de Steenrod, on peut lui associer un $GL_n(A)$-fibré différentiable en composant simplement les fonctions de transition $g_{ji}$ par la représentation $\rho$. Soit maintenant A(X) l'algèbre des fonctions[1] $C^\infty$ sur X à valeurs dans A. La méthode décrite en 3.1 permet d'associer à cette situation un projecteur dans l'algèbre des matrices sur A(X). A condition de remplacer les produits tensoriels ordinaires par des produits tensoriels complétés dans la définition de l'homologie cyclique et en utilisant la formule de Künneth, on voit alors que la caractère de Chern de E a une valeur dans la somme directe

$$\bigoplus_{r+s=2p} H^r(X) \otimes \overline{H}_s(A)$$

où $H^r(X)$ désigne la cohomologie de de Rham. On trouvera plus de détails sur ces classes caractéristiques généralisées dans la référence [ K]. Voir aussi [K2] pour les classes caractéristiques secondaires de fibrés holomorphes, algébriques ou feuilletés.

---

[1] Plus précisément le produit tensoriel topologique projectif de $C^\infty(X)$ par A.